\documentclass[a4paper,10pt]{article}
\usepackage[utf8]{inputenc}
\usepackage{amsmath,amsthm}
\usepackage{graphicx,hyperref}
\usepackage{algorithm,verbatim}
\usepackage{algpseudocode}
\newtheorem{theorem}{Theorem}[section]

\newtheorem{result}[theorem]{Result}
\newtheorem{lemma}[theorem]{Lemma}

\newtheorem{cor}[theorem]{Corollary}
\newtheorem{conj}[theorem]{Conjecture}
\newtheorem{que}[theorem]{Question}

\newtheorem{remark}[theorem]{Remark}

\title{The horizon of $2$-dichromatic oriented graphs}
\author{
J. Barát\thanks{Research supported by ERC Advanced Grant "GeoScape" and NKFIH Grant K. 131529.}\\
\small Alfr\'ed R\'enyi Institute of Mathematics\\
\small and\\
\small University of Pannonia, Department of Mathematics\\
\small 8200 Veszprém, Egyetem utca 10., Hungary\\
 and \\ 
Mátyás Czett\thanks{Supported by the ÚNKP-20-6 New National Excellence Program of the Ministry for Innovation and
Technology from the source of the National Research, Development and Innovation Fund.}\\
\small Eötvös Loránd University\\
\small 1117 Budapest, Pázmány Péter sétány 1/C, Hungary
}
\date{}

\begin{document}

\maketitle

\begin{abstract}
The dichromatic number of a directed graph is at most 2, if we can 2-color the vertices such that each monochromatic part is acyclic.
An oriented graph arises from a graph by orienting its edges in one of the two possible directions.
We study oriented graphs, which have dichromatic number more than 2.
Such a graph $D$ is $3$-dicritical if the removal of any arc of $D$ reduces the dichromatic number to 2.
We construct infinitely many $3$-dicritical oriented graphs.
Neumann-Lara found the four $7$-vertex $3$-dichromatic tournaments.
We determine the $8$-vertex $3$-dichromatic tournaments, which do not contain any of these, there are $64$ of them.
We also find all $3$-dicritical oriented graphs on $8$ vertices, there are $159$ of them.
We determine the smallest number of arcs that a $3$-dicritical oriented graph can have.
There is a unique oriented graph with $7$ vertices and $20$ arcs.
\end{abstract}

\section{Preliminaries}

The in-degree(out-degree) of a vertex in a digraph is the number of incoming(outgoing) arcs.
The degree of a vertex in a digraph is the sum of the in- and out-degrees.

A $k$-coloring of a digraph $D$ is a mapping of the vertices to $\{1,\dots ,k\}$ such that each color class is acyclic.
To be more perceptible, in the arguments with 2 colors, we usually refer to the colors as red and blue. 
In what follows, we concentrate on oriented graphs, that is directed graphs without $2$-cycles. 
We also exclude loops and parallel arcs.
 The dichromatic number $\chi_d$ of a digraph $D$ is the minimum integer $k$ such that $D$ admits a $k$-coloring.
 A digraph $D$ is $k$-dicritical if $\chi_d(D)\ge k$, but $\chi_d(D')<k$ for every proper subdigraph $D'$ of $D$.
 Let $Crit(k)$ denote the class of $k$-dicritical digraphs.
 Now $Crit(2)$ consists of all circuits (directed cycles).
 For undirected graphs, $3$-critical graphs are characterized as the odd cycles. 
 However, it is widely open to describe $Crit(3)$ for digraphs.
 Neumann-Lara described the $3$-dichromatic tournaments on 7 vertices in \cite{nl2}.
 {}From this, it is easy to find the $3$-dicritical oriented graphs on 7 vertices. 
 

The most notable conjecture regarding the dichromatic number of oriented graphs was 
raised by Neumann-Lara in \cite{n-l} and independently by \v Skrekovski \cite{circular}.

\begin{conj}\label{planar2}
Every orientation of a planar graph has dichromatic number at most $2$.
\end{conj}

For instance, the dual of the dodecahedron contains two induced paths $P_1$ and $P_2$, that together cover all vertices. 
If we color the vertices of $P_1$ red and the vertices of $P_2$ blue, then there is no monochromatic cycle. 
Therefore, any orientation of the icosahedron has dichromatic number at most 2.

It is tempting to ask whether the same method works and conclusion holds for the dual of fullerenes (and other graph classes).
The vertex-arboricity $a(G)$ of a graph $G$ is the minimum number of subsets into which $V(G)$ can be partitioned so that each subset induces a forest.
The above coloring argument works for any planar graph with arboricity at most 2.
Stein proved that a planar triangulation $T$
satisfies $a(T)\le 2$ if and only if $G^*$ is Hamiltonian \cite{stein}.
On the other hand, Kardo\v s proved that 3-connected planar cubic\footnote{$3$-regular} graphs with faces of size at most 6 are Hamiltonian.
These premises hold for fullerenes \cite{feri}.
That is, duals of fullerenes satisfy Conjecture~\ref{planar2}. 

\begin{cor}
 The dual of any fullerene has dichromatic number at most $2$.
\end{cor}

Li and Mohar \cite{li} proved that oriented planar graphs of digirth 4 are 2-dichromatic.
Harutyunyan and Mohar \cite{twores} proved there exist bounded degree digraphs with large dichromatic number and large girth.
The very natural generalisation to list-coloring has been addressed in several papers \cite{ben,har}.

 Our main target was to study the landscape of 2-dichromatic oriented graphs.
 This landscape has a boundary or horizon. 
 We consider the critical 3-dichro\-matic graphs in the next Sections.
 These are the graphs just outside of the landscape.
 
 In various related questions, the boundary of the landscape can be described by the saturated graphs of the class.
 A (di)graph $D$ is {\em saturated} in the class $\mathcal{C}$ if $D$ belongs to $\mathcal{C}$, but if we add any edge/arc to $D$, then the resulting graph $D'$ does not belong to $\mathcal{C}$.
 In our context of saturation of oriented graphs, we exclude tournaments since there are no arcs to be added.
We fix $\mathcal{C}$ as the $2$-dichromatic graphs.
 Here we notice the following interesting phenomenon.
 
 \begin{result}
 Saturated $2$-dichromatic oriented graphs do not exist. 
 \end{result}
 
 \begin{proof}
 Assume to the contrary that $D$ was $2$-dichromatic and adding any arc to $D$ results in an oriented graph that is not $2$-dichromatic.
We know that the vertices of $D$ can be split into two sets $R$ and $B$ such that each of them induces an acyclic digraph.
Let $a'$ be a missing arc from $u$ to $v$ in $D$, and let us add $a'$ to $D$ to form $D'$.
Split the vertices of $D'$ to $R$ and $B$ again.
By the assumption, one of the classes, $R$ say, induces a directed cycle $C'$, that necessarily contains the arc $a'$. 
Let $P_1$ be the directed path $C'-a'$ from $v$ to $u$.
Add now the reversed arc $a''$ from $v$ to $u$ in $D$ to form $D''$.
Split the vertices of $D''$ to $R$ and $B$ again.
By the assumption, one of the classes induces a directed cycle $C''$, that necessarily contains the arc $a''$.
We know that $u$ and $v$ were in $R$, the cycle $C''$ lies also in $R$.
Let $P_2$ be the directed path $C''-a''$ from $u$ to $v$.
Consider now the union of $P_1$ and $P_2$. 
It necessarily contains a circuit, that lies entirely in $R$, so monochromatic.
This contradicts that $D$ was $2$-dichromatic.
\end{proof}

\begin{remark}
The proof generalises to any number $k$ of colors instead of $2$.
\end{remark}


Therefore, we are left with the following

 \begin{que}
  Which oriented graphs have dichromatic number $3?$ Which are the critical ones$?$
 \end{que}
 
 Among other things, we show several infinite families of $3$-dicritical graphs and determine all of them on at most 8 vertices.
 We notice here that each infinite class collapses under taking butterfly minors.
 This concept refers to the following operation, that preserves the directed cycles of an oriented graph, and therefore seems relevant to dicoloring.
 A {\em butterfly minor} of a digraph $G$ is a digraph obtained from a subgraph of $G$ by contracting arcs, which are either the only outgoing arc of their tail or the only incoming arc of their head \cite{jrst}.
 So it feels plausible to constrain $3$-dicoloring to butterfly minor-minimal graphs rather than just critical graphs.
 However, the situation might be more subtle\footnote{Also the parity of the length of the circuits might be important. See Lemma~\ref{elso}.}.
 During a butterfly minor operation, deleting edges might result in a $2$-dicolorable graph.
 This might be later corrected by a butterfly contraction.
 We describe such an example in the penultimate paragraph of Section~\ref{crit3}.
 

 Recall that $a(G)\le 2$ implies that every orientation $D$ of $G$ satisfies $\chi_d(D)\le 2$. 
We simply color the partition classes monochromatic.
This way we cannot create a monochromatic circuit, independent of the orientations.
Kronk and Mitchem \cite{km75} proved that $a(G)\le \lceil\frac{\Delta(G)}{2}\rceil$ unless $G$ is an odd clique or a cycle.
In particular, this implies the following

\begin{cor} \label{4-reg}
For every $4$-regular graph $G$, any orientation $D$ of $G$ satisfies $\chi_d(D)\le 2$.
\end{cor}

Dirac \cite{dir} proved that any $k$-critical graph has minimum degree $k-1$.
Similarly, the minimum in- and out-degree of a $k$-dicritical graph is at least $k-1$.
Therefore, any $3$-dicritical graph satisfies $\delta^-\ge 2$ and $\delta^+\ge 2$, and need to have more than $2n$ edges by Corollary~\ref{4-reg}.
 
 The degree condition can be sometimes changed to degeneracy.
 However, Corollary~\ref{4-reg} does not generalise to 
 $4$-degenerate graphs as Lemma~\ref{elso} shows.

\section{Infinitely many critical 3-chromatic digraphs} \label{crit3}

\begin{figure}[!h]
 \begin{center}
  \includegraphics[width=0.25\linewidth]{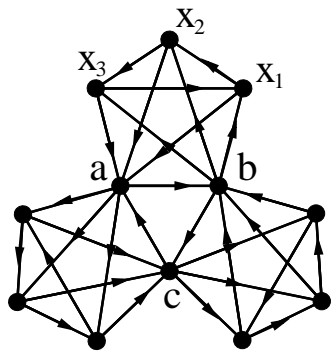}
  \caption{A $3$-dicritical oriented graph using many directed triangles.}
\label{3moons}
 \end{center}
 \end{figure}

We define a simple 5-vertex gadget $S$.
Let $x_1x_2x_3$ be a directed triangle and $ab$ an arc.
Let $b$ dominate each $x_i$ and let each $x_i$ dominate $a$.
Altogether they form an orientation of $K_5$.
We glue together 3 copies of $S$ at the 3 vertices of a directed triangle $abc$ such that each copy of $S$ --- call them $S_{ab}$, $S_{bc}$, $S_{ca}$ --- contains one of the arcs $ab$, $bc$ and $ca$ as also shown in Figure~\ref{3moons}.
The resulting oriented graph $D_1$ is $3$-dicritical as we show next.

\begin{lemma} \label{elso}
 The oriented graph $D_1$ is a member of $Crit(3)$.
 It has $12$ vertices, $9$ of them have degree $4$, and $3$ vertices have degree $8$ in the underlying graph.
 There are altogether $30$ arcs.
\end{lemma}

\begin{proof}
Suppose to the contrary that we 2-color the vertices.
 Consider the directed triangle $abc$ in the middle. 
 One of the arcs, $ab$ say, is monochromatic, red say.
 The end-vertices $a$ and $b$ are adjacent to three vertices of a directed triangle, $xyz$ say.
 If any of $x,y,z$ is red, then $a,b$ together with this red vertex form a red directed triangle.
 Otherwise $xyz$ is monochromatic blue.
 This shows that the dichromatic number of $D_1$ is at least 3.
 
 We have to show criticality.
 By symmetry, we have to consider three cases. 
 Either we delete an arc of the triangle in the middle or an arc of an outer triangle or one that connects the two types.
 The latter can also go in two directions.
 
 Suppose the deleted arc is $ab$ in the middle triangle $abc$.
 We color $a$ and $b$ red and $c$ blue.
 Now there is no monochromatic circuit containing any two of $a,b,c$ as consecutive vertices.
 We use two colors on the three vertices of any outer triangle.
Now, there is no monochromatic circuit in any of the copies of the gadget $S$.
 
Suppose next the deleted arc was $x_1x_2$ in the outer triangle of $S_{ab}$.
 We color $a$ and $b$ red and $c$ blue.
Now there is no monochromatic circuit containing $bc$ or $ca$.
  We use two colors on the three vertices of the outer triangle connected to $bc$ and $ca$.
  We color the 3 vertices of the outer triangle in $S_{ab}$ blue.
 Now, there is no  monochromatic circuit in any of the copies of the gadget $S$.
 
 Suppose last that the arc $bx_1$ is missing in $S_{ab}$.
 We color $a$, $b$ and $x_1$ red and $x_2,x_3$ blue.
 Therefore, this copy of $S$ contains no monochromatic circuit.
 We color $c$ blue and use two colors in the remaining two outer triangles.
 Therefore, there is no monochromatic circuit in the other two gadgets either.
\end{proof}

We notice the essence of the previous proof, and remark the directed triangle in the middle may be replaced by any odd circuit $O$ and the outer circuits might have any length. 
The connection between the middle and the outer circuits must remain that each arc of $O$ is contained in a gadget similar to $S$.
For each arc $ab$ of $O$, the vertices of the outer cycle $C_{ab}$ dominate $a$ and vertex $b$ dominates the vertices of $C_{ab}$.

\begin{cor} \label{vegtelen}
 There are infinitely many non-isomorphic oriented graphs in $Crit(3)$.
\end{cor}

One might say that these infinitely many graphs are very similar and argue that changing/reducing the length of the outer cycle should be allowed.
Indeed, with respect to butterfly minors, this infinite class collapses to one graph.

However, since the length of the middle circuit is an arbitrary odd number, we still have infinitely many different $3$-dicritical graphs.

We describe another construction.
Let $v_1v_2v_3$ and $v_4v_5v_6$ form two disjoint directed triangles.
We construct $D_3$ as follows.
Let all arcs go from $\{v_1,v_2,v_3\}$ towards $v_7$.
Let all arcs go from $\{v_4,v_5,v_6\}$ towards $\{v_1,v_2,v_3\}$.
Finally let all arcs go from $v_7$ towards $\{v_4,v_5,v_6\}$.
It is easy to check that $D_3$ belongs to $Crit(3)$.

We can generalise the above construction to an infinite family to further support Corollary~\ref{vegtelen}.
Instead of the two disjoint triangles let $x_1,\dots,x_k$ and $y_1,\dots,y_l$ be two circuits of any length greater than 2.
Let $v$ be a special vertex such that $v$ dominates each $x_i$ and each $y_j$ dominates $v$.
We also draw every arc of form $x_iy_j$.
The oriented graphs $O_{k,l}$ constructed this way are all $3$-dicritical.
At the same time, we again notice that this class collapses under the butterfly minor operation.
We delete an arc $vx_2$.
Now $x_1x_2$ is the only arc entering $x_2$.
Therefore, we butterfly contract the arc $x_1x_2$ and get a construction as before, but the length of the $x$-circuit is one less\footnote{We can delete the parallel arcs.}.
Notice here the following interesting fact:
The graph $O_{k,l}-vx_2$ is 2-dicolorable. 
Indeed, let $v,x_2,y_2,\dots y_l$ be blue and all other vertices red.
Now the monochromatic triples containing $v,x_2$ and $y_j$ do not form a circuit any more. 
However, after the butterfly contraction, $O_{k,l-1}$ is 3-dichromatic again.

Let us delete the matching $v_4v_1$, $v_5v_2$ and $v_6v_3$ from $D_3$ and add two arcs $v_2v_5$ and $v_3v_6$ to get the oriented graph $D_2$.
One can check that $D_2$ is $3$-dicritical.
We can generalise the construction of $D_2$ as follows.
Let $u_1\dots u_k$ and $v_1\dots v_k$ be two disjoint circuits for $k\ge 3$. 
Let $x$ be a vertex that dominates $\{v_1,\dots,v_k\}$ and each of $\{u_1,\dots,u_k\}$ dominates $x$. 
Add all arcs of form $v_iu_j$, where $i\neq j$.
Also add two arcs $u_1v_1$ and $u_2v_2$.
One can prove that these oriented graphs are also 3-dichromatic.



\section{Critical graphs on eight vertices} \label{sec8}

We wrote a computer program to find each 3-dicritical oriented graph on 8 vertices.
We represented each oriented graph as a matrix $A$ with 8 rows and 8 columns labelled from 1 to 8 such that\\
$\bullet$ if $ij$ is an arc of the digraph, then the element in row $i$ and column $j$ is +, the element in row $j$ and column $i$ is –.\\
$\bullet$ if there is no arc between two vertices $i$ and $j$, then both elements $a_{ij}$ and $a_{ji}$ are 0 (there are only $0$s in the main diagonal).

We first sought tournaments on 8 vertices, which do not contain any of the three 7-vertex 3-dicritical graphs and cannot be colored with two colors.
We started with the 6440 tournaments on 8 vertices and removed vertices one at a time in each possible way. 
We dropped the 8-vertex tournaments, which contained one of the four tournaments on 7 vertices with dichromatic number 3. 
{}From the remaining tournaments, we selected those, which contained a monochromatic circuit for any vertex coloring by two colors. See Algorithm 1. 
We found 64 tournaments that met the above criteria.

Next, we started from the pool of 64 tournaments, and we looked for graphs, which cannot be 2-dicolored, but after removing any arc become 2-dichromatic. 
Notice, we might obtain the same graph in different forms.
Therefore, we checked the digraphs that had the same in-degrees and out-degrees and examined for isomorphism. 
As a result, we obtained 159 different 3-dicritical oriented graphs on 8 vertices. 
One of them had 21 arcs, see Figure~\ref{cross3}, and all others had more arcs: 11 of them had 22 arcs, 84 had 23, 51 had 24 and 12 had 25 arcs.

In various further questions, for instance filtering oriented graphs on at least 9 vertices for $2$-dicoloring, the following characterization might be useful.
What is the list $D_1,\dots D_s$ of 8-vertex oriented graphs that satisfy the following? 
If an oriented graph $O$ contains a 3-dicritical tournament on 8 vertices as a subgraph, then it contains at least one of $D_1,\dots D_s$.
We determined this list and found that 9 oriented graphs are enough.
It is worth mentioning that the containment is distributed very unevenly.
There is a graph, which is contained in only 2 tournaments on 8 vertices.
On the other hand, there is another graph, which is contained in 34 tournaments on 8 vertices.
The 9 oriented graphs are listed here in matrix form:

\begin{verbatim}
0+000-+-    0-000++-    0+000-+-
-0000++-    +0----0+    -0000++-
000+-+--    0+0+-+--    000+-0-+
00-0++--    0+-0++--    00-0--++
00+-0+--    0++-0+--    00++0+--
+----0++    -+---0++    +-0+-0+-
--+++-0+    -0+++-0+    --+-+-0+
+++++--0    +-+++--0    ++--++-0


0000-++-    0000++--    000-+0+-
00-+0-+-    000+-+--    00+00-+-
0+00-+-+    000-+-+-    0-00-+-+
0-00+++-    0-+00+-+    +000-++-
+0+-0+--    -+-00++-    -0++0+--
-+---0++    --+--0++    0+---0++
--+-+-0+    ++-+--0+    --+-+-0+
++-++--0    +++-+--0    ++-++--0


0++000--    000-+0+-    000+-0+-
-0+000-+    00+00-+-    00-0+-+-
--0---++    0-00-++-    0+0-0-++
00+0-++-    +000-+-+    -0+0-+-+
00++0-+-    -0++0+--    +-0+0--+
00+-+0+-    0+---0++    0++-+0--
++----0+    ---++-0+    ---+++0-
+--+++-0    +++-+--0    ++---++0
\end{verbatim}

\begin{algorithm}[H] \label{2col}
\caption{Do two colors suffice?}
\begin {algorithmic}[1]
\State $t$:=the adjacency matrix of the oriented graph $G$
\For {each possible $2$-coloring of the vertices of $G$}
\State $b$:=true
\For {each color}
\State $a$:=the set of the vertices of the chosen color
\State $m$:=empty set
\For {each vertex $f$ in graph $G$}
\If {ContainsCycle($t, a, m, f$)}
\State $b$:=false
\EndIf
\EndFor
\EndFor
\If {$b$}
\State the graph is $2$-dichromatic
\State exit Algorithm
\EndIf
\EndFor
\State the graph is not $2$-dichromatic
\end{algorithmic}
\end{algorithm}
\begin{algorithm}[H]
\caption{ContainsCycle($t, a, m, f$)}
\begin {algorithmic}[1]
\State $m1$:=union($m$, $f$)
\For {each vertex $v$ in $a$}
\If {$t(f, v)$==$1$}
\If {$m$ contains $v$}
\State ContainsCycle:=true
\State exit Algorithm
\Else
\If {ContainsCycle($t, a, m1, v$)}
\State ContainsCycle:=true
\State exit Algorithm
\EndIf
\EndIf
\EndIf
\EndFor
\State ContainsCycle:=false
\end{algorithmic}
\end{algorithm}

\section{Discussion}

\begin{figure}[!h]
 \begin{center}
    \includegraphics[width=\linewidth]{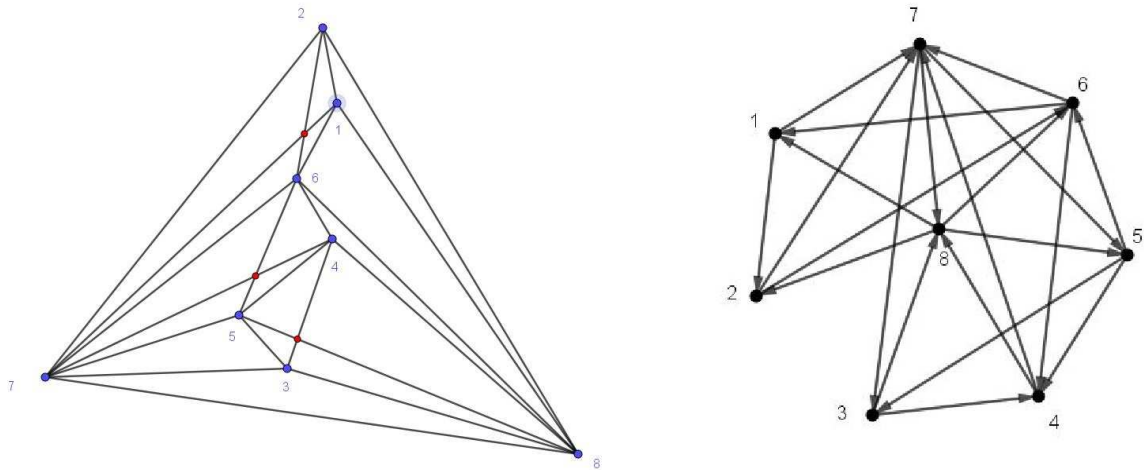} 
  \caption{The $8$-vertex $3$-critical oriented graph with $21$ arcs drawn with $3$ crossings.}
\label{cross3}
 \end{center}
 \end{figure}

It is a natural test to try strengthening Conjecture~\ref{planar2}.
One way to generalise planar graphs is to allow a few crossings in the planar embedding of $G$.
We know an oriented $3$-dicritical graph $D$, whose underlying graph $G$ has crossing number 3.
For instance, the $8$-vertex $3$-dicritical graph with 21 arcs, see Figure~\ref{cross3}.

\begin{que}
Are there $3$-dichromatic graphs with crossing number $1$ or $2$?
\end{que}

Among the $3$-critical digraphs, that we found on 7 and 8 vertices, one had 20 arcs, and three had 21 arcs. 
By the Dirac-type condition on critical digraphs, we know that the in-degree and the out-degree must be at least 2.
Therefore, if the number of vertices is at least 10, then any critical digraph must have at least 20 arcs.
By Lemma~\ref{4-reg}, we deduce that the underlying graph cannot be 4-regular either. Therefore, there must be at least 21 arcs if $n\ge 10$.
If we would like to determine the minimum number of arcs in a $3$-dicritical graph, then it remains to check the possible $9$-vertex graphs.
Again, we deduce that there must be at least 19 edges in the underlying graph.
We generated all such undirected graphs using {\em nauty} by Brendan McKay \cite{nauty} and oriented the edges to have minimum in- and out-degree 2.
We found 33700 such oriented graphs and checked its $2$-colorings to find that all of them have dichromatic number at most 2.
Similarly, we generated  
%
%
the oriented graphs with 20 arcs and 9 vertices, which have minimum in- and out-degree 2 using nauty.
%
%
%
%
%
%
There were 721603 such oriented graphs.
Again, checking the 2-colorings, we found that each of them has dichromatic number at most 2.
As a consequence, we deduce

\begin{cor}
There is only one $3$-dicritical oriented graph with $20$ arcs.
Any other member of $Crit(3)$ must have more arcs.
\end{cor}

%
%


\begin{thebibliography}{99}
\bibitem{bang}J. Bang-Jensen, T. Bellitto, T. Schweser, M. Stiebitz.
Hajós and Ore con\-struc\-tions for digraphs.
{\em The Electronic J. Combin.} {\bf 27} (1), (2020).\\
\url{https://www.combinatorics.org/ojs/index.php/eljc/article/view/v27i1p63}

\bibitem{ben}J. Bensmail, A. Harutyunyan and N.K. Le.
List coloring digraphs.
{\em J. Graph Theory} {\bf 87} (4) 492--508. (2018).


\bibitem{circular}
D. Bokal, G. Fijav\v z, M. Juvan, P.M. Kayll, B. Mohar. 
The circular chromatic number of a digraph.
{\em J. Graph Theory} {\bf 46} (3) 227--240 (2004).


\bibitem{arbo} G. Chartrand and H.V. Kronk.
The point-arboricity of planar graphs.
{\em J. London Math. Society} {\bf 44} (1)
612--616, (1969).

\bibitem{dir} G. Dirac.
The  number  of  edges  in  critical  graphs.
{\em J. für die reine und angewandte Mathematik}, {\bf 268--269} 150--164, (1974).


\bibitem{twores} A. Harutyunyan and B. Mohar.
Two results on the digraph chromatic number.
{\em Discrete Mathematics}
Volume 312, Issue 10, 1823--1826, (2012).

\bibitem{har} A. Harutyunyan and B. Mohar.
Gallai’s Theorem for List Coloring of Digraphs.
{\em SIAM J. Discrete Mathematics}
{\bf 25}(1), 170--180, (2011).


\bibitem{jrst} T. Johnson, N. Robertson, P.D. Seymour, R. Thomas. 
Directed tree-width. 
{\em J. Comb. Theory, Ser. B} {\bf 82}(1), 138--154, (2001).

\bibitem{feri} F. Kardo\v s.
A Computer-Assisted Proof of the Barnette-Goodey conjecture: Not Only Fullerene Graphs are Hamiltonian.
{\em SIAM J. Disc. Math.} {\bf 34} (1) 62--100, (2020).

\bibitem{km75} H.V. Kronk, J. Mitchem. 
Critical point-arboritic graphs.
{\em J. London Math. Soc.} {\bf 9} 459--466, (1974/75).

\bibitem{li} Z. Li and B. Mohar.
 Planar digraphs of digirth four are $2$-colourable.
 {\em SIAM J. Disc. Math.} 31(3): 2201--2205, (2018).
 
\bibitem{nauty}   B.D. McKay and A. Piperno. 
Practical graph isomorphism II.
{\em J. Symbolic Comput.}, {\bf 60} 94--112, (2014).

\bibitem{n-l}V. Neumann-Lara.
The Dichromatic Number of a Digraph.
{\em JCTB} {\bf 33}, 265--270. (1982).
\url{https://core.ac.uk/download/pdf/82562368.pdf}

\bibitem{nl2}V. Neumann-Lara.
 The $3$ and $4$-dichromatic tournaments of minimum order.
{\em Discrete Mathematics}
Volume 135, Issues 1–3, 233--243, (1994).

\bibitem{stein} S.K. Stein.
B-set and planar maps.
{\em Pacific J. Math.} {\bf 37} 217--224, (1971).

\bibitem{rs} R. Steiner.
A Note on Graphs of Dichromatic Number 2.
{\em Discrete Mathematics and Theoretical Computer Science}
{\bf 22}:4 (2021).


\end{thebibliography}
\end{document}